\documentclass{amsart}

\newtheorem{theorem}{Theorem}[section]

\theoremstyle{definition}
\newtheorem{definition}[theorem]{Definition}

\theoremstyle{remark}

\numberwithin{equation}{section}

\begin{document}

\title{Eichler Cohomology of Generalized Modular Forms of Real Weights}

\author{Wissam Raji}
\address{Department of Mathematics, American University of Beirut,
Beirut, Lebanon}


\date{}


\begin{abstract}
In this paper,  we prove the Eichler cohomology theorem of weakly
parabolic generalized modular forms of real weights on subgroups of
finite index in the full modular group. We explicitly establish the
isomorphism for large weights by constructing the map from the space
of cusp forms to the cohomology group.
\end{abstract}

\maketitle



\section{Introduction}
Eichler cohomology is the center of attention in this paper.  In
\cite{R2}, we prove Eichler cohomology theorem for generalized
modular forms of large integer weights using supplementary series.
In \cite{KLR}, we use Stokes's theorem to prove similar results on
integer weights with modifications on the multiplier system.
However, when one tries to use the classical methods to derive
Eichler cohomology theorems for real weights, obstacles with the
definition of supplementary series arise.  In \cite{KR}, we extend
the methods used in  \cite{KLR} to prove the isomorphism theorem for
real weights without any restriction on the weight.  Note that in
\cite{KR}, the isomorphism from the space of cusp forms to the
cohomology group is not as naturally defined as the corresponding
isomorphism in this paper. We use similar methods as in \cite{K4} to
derive the same results for parabolic generalized modular forms of
real weights.
\begin{definition}
A generalized modular form belonging to a subgroup $\Gamma$ of
finite index in the full modular group of real weight $k$ and
multiplier system $v$ is a function $F(z)$ satisfying:
\begin{enumerate}
\item{$F(z)$ is analytic in the upper half plane $\mathbb{H}$}
\item{ $F(z)$ satisfies a transformation law
\begin{equation}\label{eq1}
F(Mz)=v(M) (cz+d)^{k}F(z)
\end{equation}
with $\mid v \mid$ is not necessarily $1$ that depends only on the
transformation $M$ where
\begin{equation*}
M=\left(\begin{array}{lcr}
\ a  &b \\
\ c  &d \\
\end{array}\right) \in \Gamma.
\end{equation*}}
\end{enumerate}
\end{definition}
Note that the multiplier system $v:\Gamma\rightarrow \mathbb{C}$
satisfies the consistency condition
\begin{equation}\label{consis1}
v(M_1M_2)(c_3z+d_3)^k=v(M_1)v(M_2)(c_1M_2z+d_1)^k(c_2z+d_2)^k,
\end{equation}
where
\begin{equation*}
M_i=\left(\begin{array}{lcr}
\ a_i  &b_i \\
\ c_i  &d_i \\
\end{array}\right) \in \Gamma
\end{equation*}
for $i=1,2,3$ and $M_3=M_1M_2$.
\par The k-th power in (\ref{eq1}) is determined by the convention
\begin{equation}\label{conv}
w^k=\mid w\mid^ke^{ik arg w}
\end{equation}
where $-\pi \leq arg w<\pi$, for $0\neq w\in \mathbb{C}$.
\par We denote by $\{\Gamma,k,v\}$ the $\mathbb{C}$-vector space
of generalized modular forms on the group $\Gamma$ of real weight
$k$ and multiplier system $v$.
\par We shall assume that our generalized modular forms are weakly
parabolic generalized modular forms which means that $\mid v(P) \mid
=1$ for all parabolic matrices $P$.  It is important to mention that
in \cite{K4}, the multiplier system is assumed to be unitary, which
means that $\mid v(M)\mid =1 $ for all matrices $M \in \Gamma$.
\par We rewrite $(\ref{eq1})$ as
\begin{equation}\label{eq2}
F(Mz) v(M)^{-1} (cz+d)^{-k}=F(z).
\end{equation}
We introduce the stroke operator which is given by
\begin{equation}\label{eq4}
(F\mid_{v}^{k})(z)=F(Mz) v(M)^{-1} (cz+d)^{-k}.
\end{equation}
Thus (\ref{eq2}) becomes
\begin{equation}\label{eq3}
F\mid_{v}^{k}=F.
\end{equation}

It follows that
\begin{equation}\label{eq5}
F\mid_{v}^{k}M_1M_2=(F\mid_{v}^{k}M_1)\mid_{v}^{k}M_2.
\end{equation}

Let $S=\left(\begin{array}{lcr}
\ 1  &\lambda \\
\ 0  &1 \\
\end{array}\right)$, $\lambda >0$, generate the subgroup $\Gamma_{\infty}$ of translations in $\Gamma$.
 Since $F$ satisfies (\ref{eq1}), then in particular
\begin{equation*}
F(z+\lambda)=v(S)F(z)=e^{2\pi i \kappa}F(z)
\end{equation*}
with $0 \leq \kappa<1$.  Thus if $F$ is meromorphic in $\mathbb{H}$
and its poles do not accumulate at infinity, $F$ has the Fourier
expansion at $\infty$
\begin{equation}\label{eq6}
F(z)=\sum_{m=-m_0}^{\infty}a_me^{2\pi i(m+\kappa)z/\lambda}
\end{equation}
valid for $y=\Im z>y_0$.   $\Gamma$ has also $t \geq 0$ inequivalent
parabolic classes.  Each of these classes corresponds to a cyclic
subgroup of parabolic elements in $\Gamma$ leaving fixed a parabolic
cusp on the boundary of $R$, the fundamental region of $\Gamma$. We
as well denote by $\bar{R}$, the closure of the fundamental region
$R$.

\par Now let $q_1, q_2, ..., q_t$ be the inequivalent parabolic
cusps(other than infinity) on the boundary of $R$ and let $\Gamma_i$
be the cyclic subgroup of $\Gamma$ fixing $q_j$, $1\leq j\leq t$.
Suppose also that
\begin{equation*}
Q_j=\left(\begin{array}{lcr}
\ *  &* \\
\ c_j  &d_j \\
\end{array}\right)
\end{equation*}
is a generator of $\Gamma_i$; $Q_j$ is necessarily parabolic.  For
$1\leq j \leq t$; put $v(Q_j)=e^{2 \pi i \kappa_j}$, $ 0  \leq
\kappa_j <1.  $  Also $F$ has the following Fourier expansion at
$q_j$:
\begin{equation}\label{eq7}
F(z)=(z-q_j)^{-k} \sum_{m=-m_j}^{\infty}a_m(j)e^{-2\pi
i(m+\kappa_j)/\lambda_j(z-q_j)},
\end{equation}
valid for $y=\Im z >y_j$.  Here $\lambda_j$ is a positive real
number called the width of the cusp $q_j$ and defined as follows.
Let
\begin{equation*}
A_j=\left(\begin{array}{lcr}
\ 0 &-1 \\
\ 1  &-q_j \\
\end{array}\right),
\end{equation*}
so that $A_j$ has determinant 1 and $A_j(q_j)=\infty$.  Then
$\lambda_j>0$ is chosen so that
\begin{equation*}
A_j^{-1}\left(\begin{array}{lcr}
\ 1 &\lambda_j \\
\ 0  &1 \\
\end{array}\right)A_j
\end{equation*}
generates $\Gamma_j$, the stabilizer of $q_j$. We let
$C^{+}(\Gamma,k,v)$ denote the space of entire weakly parabolic
generalized modular forms of real weight $k$ and multiplier system
$v$ on $\Gamma$ which in addition to being holomorphic in $H$, it
has only terms with $m+\kappa\geq 0$ in (\ref{eq6}) and $m+\kappa_j
\geq 0$ in $(\ref{eq7})$ for all $1\leq j \leq t$. Also, let
$C^{0}(\Gamma,k,v)$ denote the subspace of weakly parabolic
generalized cusps forms which is a subspace of $C^{+}(\Gamma,k,v)$
but it has only terms with $m+\kappa> 0$ in (\ref{eq6}) and
$m+\kappa_j > 0$ in $(\ref{eq7})$ for all $1\leq j \leq t$ .
\section{Eichler Cohomology for Generalized Modular Forms of Integer Weights}
We start with definitions inherited from the classical theory of
Eichler cohomology. We now state what is known by Bol's identity
\cite{B}. With $ M=\left(\begin{array}{lcr}
\ *  &* \\
\ c  &d \\
\end{array}\right)
 \in \Gamma$,
\begin{equation}\label{eq40}
\frac{d^{k+1}}{dz^{k+1}}\left\{(cz+d)^kF(Mz)\right\}=(cz+d)^{-k-2}F^{(k+1)}(Mz)
\end{equation}
where $F^{(k+1)}$ denotes the $k+1$-th derivative of $F$. From
(\ref{eq5}), we see that if $F \in \{\Gamma,-k,v\}$ then $F^{(k+1)}
\in \{\Gamma, k+2, v\}$. The converse leads to what is called
Eichler integrals (i.e. the polynomial periods). Actually Eichler
integrals of weight $-k$ emerge when one takes the $(k+1)$-fold
integral of a modular form of weight $k+2$.  As a result,
\begin{equation}\label{eq41}
F\mid_{v}^{-k}M=F(z)+p_M(z)
\end{equation}
where $p_M(z)$ is a polynomial of degree less than or equal to $k$.
\begin{definition}
F(z) satisfying $(\ref{eq41})$ is called an Eichler integral.
Besides $\{p_M\}$ is called the system of period polynomials of $f$
(or $F$).
\end{definition}
Note that $\{p_M\mid M \in \Gamma\}$ occurring in $(\ref{eq41})$
satisfy the following consistency condition
\begin{equation}\label{eq42}
p_{M_1M_2}=p_{M_1}\mid_{v}^{-k}M_2+p_{M_2}
\end{equation}
for all $M_1$ and $M_2$ in $\Gamma$.
\begin{definition}
If $C_1=\{p_M\mid M\in \Gamma\}$ is any collection of polynomials of
degree less than or equal to $k$ such that $(\ref{eq41})$ is
satisfied, call $\{p_M\mid M\in \Gamma\}$ cocylce.
\end{definition}
\begin{definition}
If $C_2=\{p_M \mid M\in \Gamma\}$ is a set of polynomials of degree
less than or equal to $k$ such that there exists a polynomial $p$
with $p_M=p\mid_v^{-k}M-p$, then $\{p_M \mid M\in \Gamma\}$ is
called a coboundary.
\end{definition}
Observe that the coboundary satisfies $(\ref{eq42})$.
\begin{definition}
Let $P_k$ denote the complex vector space of polynomials of degree
less than $k$.  Then $H_{-k,v}(\Gamma,P_k)$ is defined by $C_1/C_2$.
\end{definition}
Parabolic cohomology plays an important role in the theory of
automorphic integrals; these are cocycles $\{p_V \mid V \in
\Gamma\}$ which satisfy the following additional condition:\\

Let $Q_0=S, Q_1,...,Q_t$ be a complete set of parabolic
representatives for $\Gamma$. Then for each $h$, $0\leq h \leq t$,
there exists a polynomial $p_h$ of degree $\leq r$ such that
$p_{Q_h}=p_h\mid Q_h -p_h$.
\begin{definition}
$\tilde{H}_{-k,v}(\Gamma,P_k)$ is a subgroup of
$H_{-k,v}(\Gamma,P_k)$defined as the space of parabolic cocycles
modulo the coboundaries.
\end{definition}

 Eichler
determined the structure of this space in terms of the space of
classical automorphic forms. His theorem for classical modular forms
\cite{E1} states that the vector space $C^+(\Gamma, k+2,v) \oplus
C^0(\Gamma,k+2, \bar{v})$ is canonically isomorphic to the first
cohomology group $H_{-k,v}(\Gamma,P_k)$.
\par In \cite{KLR},  we obtain a new result on
generalized modular forms related to Eichler cohomology.  We map the
vector space $C^+(\Gamma,k+2,v)\oplus C^0(\Gamma,k+2,\hat{v})$ into
a modified first cohomology group.  Here the action of $\Gamma$ on
$P_k$ is again by way of the slash operator in weight $-k$, but the
multiplier system is modified to $\hat{v}=vv_E$, where $v_E$ is the
multiplier system of a nontrivial entire weakly parabolic
generalized modular form of weight 0 and where $v^*$ is a unitary
multiplier system in weight $k+2$.
\section{Eichler Cohomology of Arbitrary Real Weight}
Let $k$ be an arbitrary real number and $v$ a multiplier system for
$\Gamma$ of weight $k$.  Let $\mathbb{P}$ be a vector space of
functions $g$ holomorphic in $\mathbb{H}$ which satisfy the growth
condition
\begin{equation*}
\mid g(z)\mid<K(\mid z\mid^{\rho}+y^{-\sigma}), \ \ \ \  \ y=\Im z
\end{equation*}
for some positive constant $K, \rho$ and $\sigma$. Since the weight
here is not necessarily in $\mathbb{Z}$, polynomials of fixed degree
cannot serve as the underlying space of functions in the definition
of the Eichler cohomology groups we study. Instead, we employ as the
underlying space the collection $\mathbb{P\ }$. This space was
introduced in \cite[p.612]{K4} in the context of the Eichler
cohomology theory for unitary (i.e. the usual) modular forms of
arbitrary real weight. It is worth mentioning that the
space $\mathbb{P}$ is preserved under the stroke operator.\\
\textbf{N.B.} Individual elements of $\mathbb{P}$ need not be
preserved under the stroke operator.

\begin{definition}
F is an automorphic integral of real weight -k with respect to
$\Gamma$ provided $F$ is meromorphic in $\mathbb{H}$, satisfies
\begin{equation}\label{ad1}
F\mid_{v}^{-k}M-F \in \mathbb{P} \ \ \ \mbox{for}  M\in \Gamma,
\end{equation}
and has left finite expansion at each parabolic point of $R$.
\end{definition}

As a direct consequence of the consistency condition (\ref{consis1})
for a MS $v$ in weight $k$, the slash operator satisfies
\begin{equation}  \label{add9}
F\mid_v^{-k}M_1M_2=(F\mid_v^{-k}M_1)\mid_v^{-k}M_2, \ \ \ \mbox{for}
\ M_1, M_2 \in \Gamma,
\end{equation}
where $F$ is any function defined on $\mathbb{H}$. In turn,
(\ref{add9}) implies both (\ref{consis1}) (put $F\equiv 1$) and the
"additive cocycle condition"
\begin{equation}\label{add10}
\rho_{M_1M_2}=\rho_{M_1}\mid_{v}^{-k}M_2+\rho_{M_2}, \ \ \
\mbox{for} \ M_1, M_2 \in \Gamma,
\end{equation}
where the $\rho_M$ are the periods of $F$ occurring in (\ref{ad1}). The set $%
\{\rho_M:M \in \Gamma \}$ is called a \textit{cocycle with respect to $%
\mid_{v}^{-k}$}. A \textit{coboundary} is a collection $\{\rho_M:
M\in \Gamma \}$ of elements in $\mathbb{P}$ such that
\begin{equation}  \label{add11}
\rho_M=\rho\mid_{v}^{-k}M-\rho, \ \ \ \mbox{for} \ M\in \Gamma,
\end{equation}
with $\rho$ a \textit{fixed} element of $\mathbb{P}$. Since every
coboundary
is a cocycle, we may define the \textit{Eichler cohomology group} $H_{-k,%
v}(\Gamma, \mathbb{P})$ as the quotient space:
cocycles/coboundaries.

We introduce the subspace $\tilde{H}_{-k,v}(\Gamma,\mathbb{P})$ of
"parabolic" cohomology classes in $H_{-k,v}(\Gamma, \mathbb{P})$. A
cocycle $\{\rho_M:M\in \Gamma\}$ in $\mathbb{P}$ is called
\textit{parabolic} if there exists $\rho_h$ in $\mathbb{P}$ with the
property
\begin{equation}  \label{add12}
\rho_{Q_h}=\rho_h\mid_{v}^{-k}Q_h-\rho_h, \ 1\leq h\leq t.
\end{equation}
(Recall that $Q_h\in \Gamma$ is a parabolic element such that $%
\Gamma_{q_h}=<Q_h,-I>$, where $\Gamma_{q_h}$ is the stabilizer of
the parabolic cusp $q_h$ in $\Gamma$). A coboundary is clearly a
parabolic cocycle, so we may form the quotient group: parabolic
cocycles/coboundaries.
The resulting subspace $\tilde{H}_{-k,v}(\Gamma,\mathbb{P})$ of $%
H_{-k,v}(\Gamma, \mathbb{P})$ is called the \textit{parabolic
Eichler cohomology group}.

\begin{theorem}
For any real number $k$ and for $v$ a weakly parabolic multiplier
system, we have
\begin{equation}
H_{-k,v}(\Gamma,\mathbb{P})=\tilde{H}_{-k,v}(\Gamma,\mathbb{P}).
\end{equation}
\end{theorem}
Theorem 1 is a restatement of Theorem 2 in \cite{K4} but with a
weakly parabolic multiplier system rather than a unitary one. For
the proof of Theorem 2 in \cite{K4}, Knopp mentions that the theorem
is a consequence of a result of B.A. Taylor and the proof is given
in \cite[pp.627-628]{K4}. After careful examination of Knopp's proof
for the above theorem for classical modular forms, it turns out that
the same proof can be adopted for Theorem 1 above.  This is due to
the fact that $\mid v(P)\mid=1$ for all parabolic matrices $P$ and
that the proof of Theorem 2 in \cite{K4} deals only with the
parabolic elements of the group $\Gamma$ considered.

\par The proof of one part of Theorem 4 is based
on the construction of a certain holomorphic function $\Phi$. Also
the construction of $\Phi$ involves the introduction of the "
Generalized Poincar\'e Series".
\begin{definition}
Let $w$ be a multiplier system for $\Gamma$. Let $k'$ be a positive
even integer large enough and let $\{\phi_V\}$ be a parabolic
cocycle of weight $-k$ which satisfies the additional condition that
$\phi_S=0$.  Then the Generalized Poincar\'e Series
$\Psi(\{\phi_V\},k,w;z)=\Psi(z)$ is given by
\begin{equation}\label{13}
\Psi(z)=\sum_{V \in \Theta}\phi_V(z)\bar{w}(V)(c z+d)^{-k'},
\end{equation}
where $\Theta$ is a complete set of coset representatives for
$\Gamma/\Gamma_{\infty}$ and $V=\left(\begin{array}{lcr}
\ *  &* \\
\ c  &d \\
\end{array}\right) $.
\end{definition}
Here $\mid w \mid$ is not necessarily 1. The estimation of the
cocycle is done in the same way as in \cite{K4}. The only difference
the multiplier system plays is in the convergence of the Poincar\'e
series. Recall that
\begin{equation}\label{14}
\mid w(V)\mid \leq K \mu(V)^{\alpha}
\end{equation}
where $K$ is a positive constant, $\alpha$ is another constant
depending on the modulus of the multiplier system at the generators
of $\Gamma$ and $\mu(V)=a^2+b^2+c^2+d^2$ where $a,b,c,d$ are the
entries of $V$. Thus after carefully examining the bound for the
cocycle in [\cite{K4}, Lemma 6], a similar bound for the cocycle can
be derived here by only taking the bound of the multiplier system
into consideration. Thus
\begin{equation}\label{15}
\mid \phi_V(z)\bar{w}(V)(cz+d)^{-k'}\mid\leq
K_1^*(c^2+d^2)^{e+1+\alpha}(\mid z\mid^{\eta}+y^{-\eta})
\end{equation}
where $z=x+iy$, $\eta=6e-2k$, $e=max(\rho/2,\sigma+k/2)$ with $\rho$
and $\sigma$ constants independent of the particular generator
involved in the estimation of the cocycles.  It is important to
mention that the Generalized Poincar\'e series converges absolutely
for $k'>\psi=2e+4+2\alpha$.
\begin{theorem}
For $k'>\psi$, the Generalized Poincar\'e Series (\ref{13})
converges absolutely and, in fact, $\Psi \in \mathbb{P}$.
\end{theorem}
Since by Theorem 2 the series converges absolutely for sufficiently
large $k'$, it follows that
\begin{equation}\label{15'}
\Psi\mid_v^{k} M=w(M)(\gamma z+\delta)^{k'}\Psi(z)-w(M)(\gamma
z+\delta)^{k'}g(z)\phi_M(z),
\end{equation}
for $M=\left(\begin{array}{lcr}
\ \alpha  &\beta \\
\ \gamma  &\delta \\
\end{array}\right) \in \Gamma$ where $g(z)$ is the generalized Eisenstein series
\begin{equation}\label{15''}
g(z)=\sum_{V \in \Theta}\bar{w}(V)(cz+d)^{-k'}.
\end{equation}
We now state and prove a theorem that is similar to Theorem 3 in
\cite{K4} but under the condition that $v$ is a multiplier system
associated with weakly parabolic generalized modular forms.
\begin{theorem} Let $k$ be any real number, $k>\psi>0$
and $v$ a weakly parabolic multiplier system of weight $k$. Suppose$
\{\phi_V\mid V\in \Gamma\}$ is a parabolic cocycle of weight $-k$ in
$\mathbb{P}$; that is, $\phi_V \in \mathbb{P},
\phi_{V_1V_2}=\phi_{V_1}\mid_v^{-k}V_2+\phi_{V_2}$, for all
$V_1,V_2\in \Gamma,$ and for each $j$ such that $0\leq j \leq t,$
there exists
$\phi_{Q_j}=\phi_j\mid_v^{-k}Q_j-\phi_j.$\\
Then there exists a function $\Phi,$ holomorphic in $H$, such that
\begin{equation}\label{16}
\Phi\mid_v^{-k}V-\Phi=\phi_V
\end{equation}
for all $V\in \Gamma$, and with expansions at the parabolic cusps
$q_j, 0\leq j\leq t,$ of the form
\begin{equation}\label{17}
\Phi(z)=\phi_j(z)+(z-q_j)^{-k}\sum_{m=-m_j}^{\infty}a_m(j)exp\left\{\frac{-2\pi
i(m+\kappa_j)}{\lambda_j(z-q_j)}\right\},
\end{equation}
for $1\leq j\leq t$,
\begin{equation}\label{18}
\Phi(z)=\phi_0(z)+\sum_{m=-m_0}^{\infty}a_m(0)exp\left\{\frac{2 \pi
i(m+\kappa)z}{\lambda}\right\}
\end{equation}
for $j=0$.
\end{theorem}
The construction of $\Phi$ involves Poincar\'e series defined above.
\begin{proof}
Let $\{\phi_V\}$ be a parabolic cocycle in $\mathbb{P}$ for $\Gamma$
of weight $-k$ and multiplier system $v$.  For $V\in \Gamma$, put
\begin{equation}\label{19}
\phi^*_V=\phi_V-(\phi_0\mid_v^{-k}V-\phi_0).
\end{equation}
Thus $\{\phi^*_V\}$ is a parabolic cocycle in $\mathbb{P}$.  Notice
that $\phi^*_S=\phi_0\mid S-\phi_0-(\phi_0\mid S-\phi_0)=0$.  As a
result we can form the Poincar\'e series
$\Psi(\{\phi^*_V\},k',w;z)=\Psi^*(z)$; and $\Psi^*(z)\in
\mathbb{P}$. Now let $F^*(z)=-\Psi(z)/g(z)$ then
$F^*\mid_v^{-k}M-F^*=\phi^*_M$, for $M \in \Gamma$. Define
$F(z)=F^*(z)+\phi_0(z)$, we have
\begin{equation}\label{20}
F\mid_v^{-k}M-F=\phi^*_M+\phi_0\mid M-\phi_0=\phi_M
\end{equation}
for $M \in \Gamma$. Following \cite{K4}, we see that $g(z)\neq 0$
and thus $F$ is meromorphic in $\mathbb{H}$.
\par Since $g$ has finitely many zeroes in $\bar{R}\cap
H$, $F$ has finitely many poles in $\bar{R}\cap H$.  For each $j$,
$0\leq j\leq t$, consider the function $F_j=F-\phi_j$. Then we have
\begin{equation*}
F_j\mid Q_j-F_j=(F\mid Q_j-F)-(\phi_j\mid
Q_j-\phi_j)=\phi_{Q_j}-\phi_{Q_j}=0;
\end{equation*}
and thus it follows that
\begin{equation*}
F_j(z)=(z-q_j)^{-k}\sum_{m=-\infty}^{\infty}a_m(j)exp\left\{\frac{-2\pi
i(m+\kappa_j)}{\lambda_j(z-q_j)}\right\},
\end{equation*}
for $1\leq j\leq t$,
\begin{equation*}
F_0(z)=\sum_{m=-\infty}^{\infty}a_m(0)exp\left\{\frac{2 \pi
i(m+\kappa)z}{\lambda}\right\}
\end{equation*}
for $j=0$.\\
Also by the same argument as in \cite{K4}, $F(z)$ has expansions of
the form (\ref{17}), (\ref{18}).
\par Since $F$ may have poles in $H$ we need to modify it to
obtain the function $\Phi$ of Theorem 3.  It follows from \cite{P1}
that there exists $f\in \{\Gamma, k,v_u\}$ with poles at given
principal parts at finitely many points of $\bar{R}\cap H$ and is
otherwise holomorphic in $R$ with possible exceptions at the cusps.
In \cite{KM3}, Knopp and Mason proved that given a multiplier system
of a generalized modular form $v$, then $v$ can be written as
$v=v_Ev_u$ where $v_E$ is the multiplier system of an entire
generalized modular form $E$ of weight $0$ and $v_u$ is unitary.
Note that $E(z)$ is an entire generalized modular form of dimension
0 that has no zeroes or poles \cite{KM3}.  The existence of such a
function E(z) can be found in \cite{K1}. As a result, $fE$ is a
generalized modular form of weight $-k$ and multiplier system $v$
that has its poles and zeroes as $f$. We form now $\Phi=F-fE$.
Notice that by results from Petersson in \cite{P1}, there exists
such an  $f\in \{\Gamma,-k,v\}$ which have poles whose principal
parts agree with those of $F$ in $\bar{R}\cap H$. Since
$fE\mid_v^{-k}M=fE$ for $M\in \Gamma$, we still have
$\Phi\mid_v^{-k}M-\Phi=\phi_M$. Since $fE$ has a pole of finite
order at each parabolic cusps, then $\Phi$ has the expansions
$(\ref{17})$ and $(\ref{18})$ at those cusps. Finally, $\Phi$ is
holomorphic in $\bar{R}\cap H$.
\end{proof}
\par We now state the main theorem in this paper.  Notice that the result here
follows from Theorem A* in \cite{KR}.  However, in the present
paper, we are explicitly constructing  the isomorphism which seems
to arise naturally and which actually fails for small weights due to
convergence issues of the generalized Poincar\'e series.
\begin{theorem}\textbf{The Main Theorem}
Let $\psi$ be an integer computable depending on $k$ and $v$. If
$k\leq -2$ or $k > \psi>0 $ with $v$ a multiplier system of weight
$k$, then
$H_{-k,v}(\Gamma,\mathbb{P})=\tilde{H}_{-k,v}(\Gamma,\mathbb{P})$
are isomorphic to $C^0(\Gamma, k+2,\bar{v})$.
\end{theorem}

\textbf{Case 1} $k\leq -2$ \\ By Mittag-Leffler theorem for
automorphic forms, there exist $G'\in \{\Gamma,-k,v_u\}$ such that
$G'$ is holomorphic in $\mathbb{H}$ and G' has a principal part at
each cusp that agrees precisely with the principal part at the cusp
of the expansion of $\Phi$. As in the proof of Theorem 3 above and
due to the result of Knopp and Mason in \cite{KM3}, we can write
$v=v_uv_E$ where $E$ is again an entire generalized modular form of
weight 0. Now let $G=G'E$, thus $G\in \{\Gamma,-k,v\}$, holomorphic
in $\mathbb{H}$ and has a principal part at each cusp as that of
$G'$. So put $\Phi^*=\Phi-G$, then $\Phi^*$ is holomorphic in
$\mathbb{H}$ and has expansion (\ref{17}) at each cusp $q_j$ in
which no negative power in the local parameter appears and as a
result, $\mid\Phi^*(z)\mid <K(\mid y\mid^{\rho}+ y^{-\sigma} )$ for
constants $K,\rho, \sigma$ and for $z=x+iy\in \mathbb{H}\cap
\bar{R}$. Also
\begin{equation*}
\Phi^*\mid_v^{-k}-\Phi^*=\phi_M, \ \ \ M\in \Gamma.
\end{equation*}
Notice also as in \cite{K4}, if we get that $\mid
\Phi^*(z)\mid<K(y^{\alpha}+y^{-\beta})$ for $z \in \bigcup_{V\in
\Gamma/\Gamma_{\infty}}V(\bar{R})\cap \mathbb{H}$ then $\Phi^*\in
\mathbb{P}$ and hence $H_{-k,v}(\Gamma,\mathbb{P})=0$.  To prove
this, we follow the same proof as in \cite[p.622]{K4} by
constructing $f(z)=y^{-k/2}\mid \Phi^*(z)\mid$, $y=\Im z>0$.  Then
for $M=\left(\begin{array}{lcr}
\ *  &* \\
\ c  &d \\
\end{array}\right)\in \Gamma,$
\begin{equation*}
f(Mz)=\mid v(M)\mid f(z)+y^{-k/2}\mid v(M)\mid \mid \phi_M(z)\mid.
\end{equation*}
The only difference here from that of \cite{K4} is that $\mid
v(M)\mid$ is not necessarily equal to 1 for non-parabolic elements M
of $\Gamma$. As a result, we get the same inequality with a change
in the estimate as in \cite[p.623]{K4} and after carefully checking
that Lemma 8 is still valid for our purposes.
\begin{equation*}
\mid\Phi^*(z)\mid\leq y^{k/2}\left(\frac{1}{my}+y\right)^{-k/2}\mid
v(V)\mid
\left\{K_R\left(\frac{1}{my}+y\right)^{\rho}+y^{-\sigma}+K_1(y^{\rho_0}+y^{-\sigma_0})\right\}
\end{equation*}
where $K_R, K_1, \sigma_0, \rho_0$ and $m$ are positive constants.\\
\par \textbf{Case 2} $k>\psi$
\par
The proof of this case follows exactly as the proof of Theorem 1 in
\cite{K4} for the case $r>0$ with two differences that will be
mentioned below.

1. The first difference is the weakly parabolic multiplier system
that is considered here. It turned out that the multiplier system
will not impact the proof if we make sure that the theorems used by
Knopp continue to work here. Specifically, in the proof of Theorem 1
in \cite{K4}, Knopp used three theorems from Niebur's \cite{N2, N1}.
Thus showing that those theorems will also work for the case of
weakly parabolic multiplier systems will be sufficient to deduce the
proof of Theorem 4 for the case $k>\psi$ in the present paper.

2. The other difference is the existence of of a function
$f\in\{\Gamma,r,v\}$ in \cite[p. 623]{K4} (the case $r>0$) which has
poles of prescribed principal parts at each of the cusps
$q_1,...,q_t$ and is holomorphic in $\mathbb{H}$.  However, the
existence of this function is also guaranteed for weakly parabolic
multiplier by multiplication of $f$ with an entire generalized
modular form $E$ of weight 0. This method is used in the present
paper in the proof of Theorem 3 and also repeated in the proof of
this theorem for the case $k\leq -2$.
\par We now list the extensions of the three theorems of Niebur that were mentioned
above in the context of weakly parabolic multiplier systems. The
Generalization of Theorem $N_1$ in the context of weakly parabolic
generalized modular forms is proved as Theorem 1 in\cite{R3}. In
fact, in \cite{R3}, we generalize two theorems from \cite{N2, N1}
(mentioned in \cite{K4} as Theorems $N_1
N_2$) to the case of weakly parabolic generalized modular forms.\\
\textbf{Generalization of Theorem $N_1$.} Let $k >\psi>0$ where
$\psi$ is an integer computable depending on $k$ and $v$. Let $m_0$
be a nonnegative integer and $a_{-1},...,a_{-m_0}$ complex numbers.
Then there exists $F$, an automorphic integral of weight $-k$ with
respect to $\Gamma$, which is holomorphic in $\mathbb{H}$ and at
finite cusps $q_1,...,q_t$ and which has an expansion (\ref{17})
with principal part
\begin{equation}\label{30}
a_{-m_0}e^{2\pi i(-m_0+\kappa)z/\lambda}+...+a_{-1}e^{2 \pi i
(-1+\kappa)z/\lambda}
\end{equation}
at $q=\infty$.  The function $F$ has the transformation properties
\begin{equation}\label{31}
F\mid_{v}^{-k}V=F, \ \ \ V\in \Gamma_{\infty}
\end{equation}
and
\begin{equation}\label{32}
\overline{F\mid_{v}^{-k}-F}=\int_{V^{-1}\infty}^{i\infty}G(\tau)(\tau-\bar{z})^{k}d\tau,
\ \ \ V\in \Gamma-\Gamma_{\infty},
\end{equation}
where $G \in C^0(\Gamma,k+2, \bar{v})$ is determined by
$a_{-m_0},...,a_{-1}$ and the path of integration is a vertical
line.
Furthermore, $F\mid V-F$ is a parabolic cocycle in $\mathbb{P}$.\\

The following theorem is a generalization of a theorem from
\cite{N2} that was mentioned in \cite{K4} as Theorem $N_2$. It is
actually the converse of Generalization of Theorem $N_1$.
\\
\textbf{Generalization of Theorem $N_2$.}Let $\psi$ be an integer
computable depending on $k$ and $v$. Let $k>\psi$ and let $G \in
C^{0}(\Gamma,k+2,\bar{v})$ then there exists an automorphic integral
$F$ satisfying $(\ref{31})$ and $(\ref{32})$ such that $F$ is
holomorphic in $\mathbb{H}$ and at $q_1,...,q_t$.\\ \\
The following theorem is again a generalization of Theorem $N_3$
from \cite{K4} but for weakly parabolic multiplier systems. It
follows as well from the main result of \cite{R1} where relations
between Fourier coefficients of generalized modular forms were
determined using the circle method.
\\
\textbf{Generalization of Theorem $N_3$.} Let $\psi$ be an integer
computable depending on $k$ and $v$. For $k>\psi$ If there exists
$\tilde{F}\in \{\Gamma,-k,v\}$ which is holomorphic in $\mathbb{H}$
and at $q_1,...,q_t$, and which has principal part (\ref{30}) at
$q_0=\infty$, then the function $F$ of Theorem A is in
$\{\Gamma,-k,v\}$ and in fact $F=\tilde{F}$.  In this case the cusp
form $G$ of (\ref{32}) is $\equiv 0$.

\bibliographystyle{amsplain}

\begin{thebibliography}{10}

\bibitem {B} G. Bol, \textit{Invarianten linearer Differentialgleichungen}, Abh. Math. Sem. Univ. Hamburg 16(1949), nos 3-4, 1-28.\\
\bibitem {E1} M. Eichler, Grenzkreisgruppen und kettenbruchartige
Algorithmen, Acta Arith. 11 (1965) 169-180.\\
\bibitem {H1} S. Husseini and M. Knopp, \textit{Eichler Cohomology of Automorphic
Forms}. Illinois J. Math 15 (1971) 565-577.\\
\bibitem {K4} M. Knopp, \textit{Some new results on the Eichler Cohomology of automorphic
forms}, Bull. Amer. Math. Soc. 80 (1974), 607-632.\\
\bibitem {K1} M. Knopp and G. Mason, \textit{Generalized modular forms}. J. Number Theory 99 (2003), 1-18.\\
\bibitem {KLR} M. Knopp, J. Lehner and W. Raji, \textit{Eichler cohomology for generalized modular forms}.
International Journal of Number Theory, Vol. 5 No. 6 (2009), 1049-1059.\\
\bibitem {KR} M. Knopp and W. Raji, \textit{Eichler cohomology for generalized modular forms II}. International Journal of Number Theory, Vol 6 No. 5 (2010) 1083-1090.\\
\bibitem {KM3} M. Knopp and G. Mason, \textit{Parabolic generalized modular forms and their characters}. International Journal of Number Theory V.5 No. 5 (2009), 845-857.\\
\bibitem {N2} D. Niebur, \textit{Automorphic Integrals of arbitrary positive dimension and Poincare series}, Doctoral Dissertation, Madison, Wis., 1968.\\
\bibitem {N1} D. Niebur, \textit{Construction of Automorphic forms and
integrals}.  Trans. of the AMS. 191 (1974) 373-385.\\
\bibitem {P1} H. Petersson, \textit{Konstruktion der Modulformen und zu gewissen Grenzkreisgruppen gehorigen automorphen
Formen von positiver reeller Dimension und die vollstandige
Bestimmung threr Fourier koeffizienten,} S.B. Heidelberger Akad.
Wiss Math.-Nat. Kl. 1950, 417-494. MR 12, 806.\\
\bibitem {R1} W. Raji, \textit{Fourier coefficients of generalized modular forms of negative weights},
International journal of number theory Vol.5, 2 (2009), 153-160.\\
\bibitem {R2} W. Raji, \textit{Eichler cohomology theorem of generalized modular
forms}, International Journal of Number Theory, Vol 7 No. 4 (2011) 1103-1113.\\
\bibitem {R3} W. Raji, \textit{Construction of Generalized Modular Integrals}, Functiones et Approximatio
Commentarii Mathematici Vol 41 No.2 (2009) 105-112.\\


\end{thebibliography}

\end{document}